\documentclass[10pt]{article}

\usepackage{amsthm}

\usepackage[latin1]{inputenc}
\usepackage[all]{xy}
\usepackage{amssymb}
\usepackage{pstricks}
\usepackage{amsmath}
\usepackage{array}
\usepackage{indentfirst}
\usepackage{hyperref}
\usepackage{graphicx}
\usepackage{amsxtra}
\usepackage{amscd}
\usepackage{amsmath,amssymb}

\usepackage{xcolor}
\usepackage[all]{xy}

\newcounter{nmdthmcnt}
\newenvironment{namedthm}[2][]{\addtocounter{nmdthmcnt}{1}%
\theoremstyle{plain}\newtheorem*{nmdthm\roman{nmdthmcnt}}{Theorem #2}%
\begin{nmdthm\roman{nmdthmcnt}}[#1]}{\end{nmdthm\roman{nmdthmcnt}}}

\usepackage[latin1]{inputenc}
\usepackage[all]{xy}
\usepackage{amssymb}
\usepackage{pstricks}
\usepackage{amsmath}
\usepackage{array}
\usepackage{indentfirst}
\usepackage{hyperref}
\usepackage{graphicx}

\newtheorem{thm}{Theorem}[section]

\newtheorem{lem}[thm]{Lemma}
\newtheorem{cor}[thm]{Corollary}
\newtheorem{prop}[thm]{Proposition}

\def\Ext{{\mathrm Ext}}
\def\ra{\longrightarrow}

\def\Hom{{\mathrm Hom}}

\numberwithin{equation}{section}

\def\GL{\mathbb {GL}}

\def\Z{\mathbb Z}
\def\F{\mathbb F}

\usepackage[all]{xy}
\usepackage{amssymb}
\usepackage{pstricks}
\usepackage{amsmath}
\usepackage{array}
\usepackage{indentfirst}
\usepackage{hyperref}
\usepackage{graphicx}

\def\F2{{\mathbb F}_2}
\def\U{\mathcal U}

\def\Z{\mathbb Z}

\def\Ker{{\rm Ker}}

\def\Nil{{\mathcal N}il}
\def\calF{{\mathcal F}}
\def\calP{{\mathcal P}}

\def\A{{\mathcal A}_2}

\def\U{{\mathcal U}}

\def\ra{\rightarrow}

\def\GL{\mathbb {GL}}
\def\Hom{\mathrm Hom}
\def\Ext{{\mathrm Ext}}
\def\coker{\mathrm coker}

\def\Hom{\mathrm Hom}

\begin{document}

\title{Some finiteness results in the category $\U$ }
\author{Nguyen The Cuong and Lionel Schwartz} 
\date {UMR 7539 CNRS,
 Universit\'e Paris 13\\
 LIA CNRS Formath Vietnam\\}

\maketitle

\begin{abstract} This note investigate some finiteness properties of the category $\U$ of unstable modules. One shows finiteness properties for the injective resolution of finitely generated unstable 
modules. One also shows  a stabilization result under Frobenius twist  for $\Ext$-groups.
\end{abstract}

\begin{section}{Introduction}

This note investigates some finiteness properties of the category $\U$ of unstable modules. For simplicity, only  the case of
the prime $2$ will be considered. The main finiteness property of $\U$ is to be locally noetherian \cite{MP67}, see also \cite{LZ86}. This means that if $M$ is a 
finitely generated unstable module any sub-module is also finitely generated.
Injective objects of the category are described in \cite{LS89}. It is natural and useful to ask for some control 
on injective resolutions. Some of the results of the note could have been written earlier. They are implicit in some places, however it is worth to present them all at the same time and
explicitly. The first result of the paper is:

\begin{namedthm}{\ref{uf}}
Let $M$ be an unstable finitely generated module. It has an injective resolution $\mathcal I^\bullet$ such that
any $\mathcal I^k$ is a finite direct sum of indecomposable injective unstable modules.
\end{namedthm}

This is to be compared with an analogous theorem in the category $\calF$ of functors from finite  dimensional 
$\F2$-vector spaces to $\F2$-vector spaces \cite{S94, FLS94, ER03}:

\begin{thm}{\label{pf}} Let $F$ be a polynomial functor taking finite dimensional values. There exists  an injective resolution $\mathcal I^\bullet$ of $M$  such that
any $\mathcal I^k$ is a finite direct sum of indecomposable injective functors.
\end{thm}

In the sequel such injective resolutions (for both categories) will be called of finite type.  The proof of \ref{uf} uses \ref{pf} but is not a direct consequence.

The link between \ref{uf} and \ref{pf} is given by the functor $f \colon \U  \ra \calF$ \cite{HLS93}.  This functor induces an equivalence $\U \ra \U/\Nil \cong \calF_\omega$ of the quotient category of $\U$
by the full subcategory $\Nil$  of nilpotent modules to the category of analytic functors. It has a right adjoint $m$, the composition $\ell = m \circ f$ 
is the localization functor away from $\Nil$. The natural
map $M \ra \ell(M)$ is initial for morphisms $M \ra L$, $L$ being $\Nil$-closed {\it i.e.} such that $\Ext_\U^i(N,L)=\{0\}$ for $i=0,1$, any $N \in \Nil$ \cite{Gab62}.

The functor $\ell$  is left exact and admits right derivatives $\ell^i$. In some interesting cases described later  the functors $\ell^i$
are computed by MacLane homology. Here is a general result:

\begin{namedthm}{\ref{lf}}
If $M$ is a finitely generated unstable module the unstable modules $\ell^i(M)$ are finitely generated unstable modules.
\end{namedthm}

This is  related to \ref{uf} but not (at least directly) equivalent.

Recall the Frobenius twist functor $\Phi$ on the category $\U$ \cite{S94}, this functor is also called "double functor", but the terminology "Frobenius functor" looks better \cite{NDHH10}. The stabilization result which follows is a direct corollary of \ref{uf}:

\begin{namedthm}{\ref{stab}} Let $M$ and $N$ be two finitely generated unstable modules. Consider the direct system induced by 
the maps $\Phi^{k+1} M \ra  \Phi^{k} M$: 
$$
\ldots \ra \Ext_\U^i (\Phi^kM,N) \ra \Ext_\U^i (\Phi^{k+1}M,N)\ra \ldots
$$
For $k$ large enough this map is an isomorphism, and the terms are isomorphic to $\Ext^i_\calF(f(M),f(N))$.
\end{namedthm}
A similar stabilization result holds 
in the category of strict polynomial functors for the map induced by the Frobenius twist \cite{FFSS99}, $F$ and $G$ denoting strict polynomial functors of the same degree, and $F^{(k)}$ being the $k$-th Frobenius twist:
the system
$$
\ldots \ra \Ext_\calP^i (F^{(k)},G^{(k)}) \ra \Ext_\calP^i (F^{(k+1)},G^{(k+1)})\ra \ldots
$$
stabilizes for $k$ large.  But in  general  the colimit is not isomorphic
to $\Ext^i_\calF(\mathcal O(F), \mathcal O(G))$, $\mathcal O \colon \calP \ra \calF$ is the forgetful functor. However, in interesting cases it is isomorphic, it is true in particular  when $F$ and $G$ are the canonical liftings of the simple objects in $\calF$ in the category
$\calP$. Another difference is that on the right hand side of the $\Ext$-group the module does come with a Frobenius twist: one has just $N$ instead of $\Phi^k(N)$.
Again in some interesting cases it is possible to replace $N$ by $\Phi^k(N)$. The main reason to keep $N$ is that  the colimit considered in the corollary  is more natural for topological applications.

It is also worth to note that it is  not clear to the authors how to prove the result using  projective resolutions. The category $\U$ is locally noetherian. Thus,  given a finitely generated module $M$, it is possible to construct a projective resolution of $M$ such that each term of the resolution is a finite direct sum of indecomposable projective modules. 
However when replacing $M$ by $\Phi^k(M)$, if the finite generation is preserved,  one has no  control on the  size  of a given term in the resolution. 

\vskip 2mm
In particular if $\Ext^i_\calF(f(M),f(N)) \cong \{0\}$ for $k$ large enough $\Ext_\U^i (\Phi^{k+1}M,N) \cong\{0\}$. This case applies for $i=1$ in  the
following case. Let $\U_n$, $n \in \mathbb N$,  be the Krull filtration on $\U$ (see the next section).
Let $M$ be a finitely generated unstable module, assume $M \in \U_n \setminus \U_{n-1}$. Let $R$ be the smallest sub-module in 
$\U_n$ such that $M/R \in  \U_{n-1}$. The following is a corollary of  \ref{stab}:

\begin{prop}\label{EMt}
For $k$ large enough the $\F2$-vector space $\Ext^1_\U(\Phi^k (R\otimes R), M)$ is trivial. 
\end{prop}
This is applied in \cite{CGS14}. This statement can be proved directly using only Steenrod operations, however it is much more tedious. 

\end{section}
\begin{section}{Recollections: the nilpotent filtration,  the Krull filtration and functors}

This section recalls briefly facts about the category $\U$. One refers mostly to \cite{S94} and \cite{K95} for all of this material.
The  subcategory ${\mathcal N}il_s$, $s \geq 0$, of $\U$ is the smallest thick subcategory stable under
colimits and containing all $s$-suspensions of unstable modules.
$$
\U = {\mathcal N}il_0 \supset {\mathcal N}il_1 \supset {\mathcal N}il_2  \supset \ldots \supset
{\mathcal N}il_s \supset \ldots
$$

\begin{prop} Any unstable module  $M$ has a convergent decreasing filtration
 $\{M_s\}_{s\geq 0}$ with $M_s/M_{s+1} \cong \Sigma^s R_s(M)$ 
 where $R_s(M)$ is a reduced unstable module, {\it i.e.}  which does not contain a non trivial
suspension.
\end{prop}

An unstable module is nilpotent if it belongs to $\Nil_1$.
The following is proved in \cite{S94} Lemma 6.1.4, see also \cite{K95}:

\begin{prop} \label{KS} Let $M$ be a finitely generated unstable module. Then the $R_i(M)$ are finitely generated
and trivial if $i$ is large enough.
\end{prop}

The category of unstable modules
  $\U$, as any abelian category,   has a   Krull filtration. There are  thick subcategories
which are stable under colimits 
$$
\U _0 \subset \U _1 \subset \U _2  \subset \ldots \subset \U
$$
defined as follows.

The category
$\U_0$ is the  largest thick sub-category generated by simple objects and stable under colimits. It
is the subcategory of locally finite modules, $M \in \U$  is locally finite if  the span over $\mathcal A_2$
of any $x \in M$ is finite.
 Having defined by induction $\U_n$ one defines $\U_{n+1}$ as follows. One introduces the quotient
category $\U/\U_n$ whose objects are the same of those of $\U$ but where morphisms in $\U$ that have
kernel and cokernel in $\U_n$ are formally inverted.
 Then  $(\U/\U_n)_0$ is defined  as above and $\U_{n+1}$ is the pre-image of this subcategory in
$\U$ {\it via} the canonical projection functor, see \cite{Gab62} for details. 
 One has

\begin{thm} Let $M \in \U$ and $K_n(M)$ be the largest sub-object of $M$ that is in $\U_n$, then $$M
= \cup_n K_n(M{})$$
\end{thm}

Let us give some examples: 
\begin{itemize}
\item
$\Sigma^kF(n) \in \U_n \setminus\U_{n-1}$, the unstable modules $F(n)$ are the canonical generators
of $\U$, generated in degree $n$ by $\iota_n$ and $\F2$-basis 
$Sq^I\iota_n$, $I$ an admissible multi-index of excess less than $n$;

\item  the $n$-th exterior power $\Lambda^n(F(1)), F(1)^{\otimes n} \in \U_n\setminus \U_{n-1}$...
 \end{itemize}

There is a characterization of the Krull filtration in terms of the functor $\bar T$ introduced by Jean
Lannes. 
The functor $T_V$, $V$ being an elementary abelian $p$-group, is  left adjoint to $M \mapsto H^*BV
\otimes M$. If $V=\F2$ $T_V$ is denoted  by $T$. As $H^*B\Z/2$ splits up, in $\U$, as  $\F2 \oplus \tilde
H^*B\Z/2$ the functor $T$ is naturally equivalent to $Id \oplus \bar T$. Below are the main
properties of  $T_V$:

\begin{thm}[\cite{La92,S94}] The functor $T_V$ commutes with colimits
(as a left adjoint). It  is exact. Moreover there is a canonical isomorphism
$$
T_V(M \otimes N) \cong T_V(M) \otimes T_V(N)
$$
\end{thm}

If $M = \Sigma \F2$ it writes as
$
T_V(\Sigma M) \cong \Sigma T_V(M)
$.

Below is the caracterisation of the Krull filtration alluded to above:
 \begin{thm}\label{krull} The following two conditions are equivalent:
 \begin{itemize}
 
 \item$M \in \U_n$ ,
 \item $\bar T^{n+1} (M)=\{0\}$.
 \end{itemize}
\end{thm}

A very nice proof of this result is in \cite{K13}, as well as in Nguyen The Cuong's thesis \cite{NTC14}.

\begin{cor}\label{krull-t}

If $M \in \U_m$ and $N \in \U_n$ then $ M \otimes N \in \U_{m+n}$.

\end{cor}

Let $\calF$ be the category of functors from finite dimensional $\F2$-vector spaces to all vector
spaces. Define a functor, \cite{HLS93},  $f \colon  \U \ra \calF$  by 
$$f(M)(V)=\Hom_{\U}(M,H^*(BV))^*=T_V(M)^0$$

Let $\calF_n$ be the sub-category of polynomial functors of degree less than $n$. It  is defined as
follows. Let $F \in \calF$, let $\Delta(F) \in \calF$
be defined by
$$
\Delta(F)(V)=\Ker(F(V \oplus \F2) \ra F(V))
$$
Then by definition $F \in \calF_n$ if and only if $\Delta^{n+1}(F)=0$. 
As an  example $V \mapsto V^{\otimes n}$ is in $\calF_n$. The following holds for any $M \in \U$
$$
\Delta(f(M)) \cong f ( \bar T(M))  
$$

Thus,  the  diagram below commutes:

$$
\xymatrixcolsep{2cm}\xymatrix{
\U_0 \ar@{>}[d]^{f}\ar@{^{(}->}[r]& \ldots\ar@{^{(}->}[r] & \U_{n-1} \ar@{>}[d]^{f} \ar@{^{(}->}[r]	& \U_n   \ar@{>}[d]^{f}   \ar@{^{(}->}[r]	 &\U \ar@{>}[d]^{f} \\
\calF_0\ar@{^{(}->}[r]            & \ldots\ar@{^{(}->}[r] &	\calF_{n-1} \ar@{^{(}->}[r]	            &\calF_n     \ar@{^{(}->}[r]	           & \calF\\	
}$$

An injective unstable module always splits up as the direct sum of a reduced module and of a nilpotent one, moreover (by definition) there are no non trivial maps from 
a nilpotent module to a reduced one (\cite{S94} chapter 2 and 3).  Thus, any injective resolution $\mathcal I^\bullet$ in the  $\U$ of an unstable module $M$ has the following properties.

\begin{prop}
 For any $k$  $\mathcal I^k$ decomposes as a direct sum  $\mathcal R^k \oplus \mathcal N^k$, the first module being reduced and the second one nilpotent. The differential 
 $\partial_k \colon \mathcal I^k \longrightarrow \mathcal I^{k+1}$ writes as 
$\left(\begin{smallmatrix} \partial^k_r & 0\cr \rho^k  &\partial_n^k\cr \end{smallmatrix}\right) \colon	\mathcal R^k \oplus \mathcal N^k  \longrightarrow
\mathcal R^{k+1} \oplus \mathcal N^{k+1} $.
The $k$-th cohomology module of the quotient complex $\mathcal R^\bullet$ is (by definition) the $k$-th derived functor $\ell^k(M)$ of the localization functor away from $\Nil$ applied to $M$. 
\end{prop}

The functor $f$ has a right adjoint $m$, the composition $\ell = m \circ f$ 
is the localization functor away from $\Nil$. The natural
map $M \ra \ell(M)$ is initial for $M \ra L$, $L$ being $\Nil$-closed {\it i.e.} such that $\Ext_\U^i(N,L)=\{0\}$ for $i=0,1$, any $N \in \Nil$ \cite{Gab62}. In particular the localization of a nilpotent unstable module is trivial.
It follows that:
 $$
 H^{k+1}(\mathcal N^\bullet) \cong \ell^k(M)
 $$
and the modules $\ell^k(M)$, $k\geq 1$, are nilpotent (this can be seen directly). Moreover as 
the functor $T$  preserves reduced modules and nilpotent modules it commutes with $\ell$ so that  $T(\ell(M)) \cong \ell(T(M))$ and more generally $T(\ell^i(M)) \cong \ell^i(T(M))$.

\begin{prop} \label{krull-loc}
If $M \in \U_n$, then $\ell^i(M) \in \U_{n-1}$, if $i>0$.
\end{prop}

Indeed, it is enough to show that $\bar T^n(\ell^i(M))$ is trivial if $i>0$, but $\bar T^n(\ell^i(M)) \cong \ell^i(\bar T^n(M))$ and the result follows from \ref{krull}.

As the tensor product of reduced injective unstable modules is still injective one has:
 
 \begin{cor}[K\"{u}nneth formula]\label{kunneth} 
 Let $M$ and $N$ be two unstable modules then:
 $$
 \ell^k(M \otimes N) \cong \oplus_{i+j=k} \, \, \ell^i(M) \otimes \ell^j(N)
 $$
 \end{cor}

One keeps the notation introduced above. Assume the resolution to be minimal, that is $\mathcal I^{k+1}$ is the injective hull of $\coker( \partial^{k-1})$.
On the other hand  $\mathcal R^{k+1}$ is the injective hull of the quotient $\mathcal R'^k$ of $\coker( \partial_r^{k-1})$ by its largest nilpotent submodule $N(\coker( \partial_r^{k-1}))$.
 The results which follow are standard homological algebra.

  \begin{prop}\label{4-terms} 
  Let $N'^k$ be the largest nilpotent submodule of $\coker( \partial^{k-1})$.   One has the following exact sequence:
   $$
  \{0\} \ra \ell^k(M) \ra \coker (\partial^{k+1}_n) \ra N'^{k+1} \ra \ell^{k+1}(M) \ra \{0\}
  $$  If the resolution is minimal the unstable module $\mathcal R^{k+1}$  is the injective hull of $R'^k$, 
$\mathcal N^{k+1}$ is the injective hull of $N'^k$.
  \end{prop}
  
  If one is given an injective resolution $\mathcal I^\bullet$ of an unstable module $M$, $f(\mathcal I^\bullet)$ is  an injective resolution of $f(M)$. If moreover one assumes
  $M$ to be finitely generated then $f(M)$ is a finite functor, using Kuhn's terminology: it is polynomial and takes finite dimensional values  \cite{K94-1}. Then \ref{pf}, implies that a minimal
  resolution of $f(M)$ is of finite type. This implies easily that the reduced part of a minimal resolution of $M$ is of finite type and in fact:
  
  \begin{prop}\label{puf} Let $M$ be a finitely generated unstable module, and $\mathcal I^\bullet$ a minimal injective resolution. Then $f(\mathcal I^\bullet)$ is a minimal injective resolution of $f(M)$.
  \end{prop}
  
 This follows from the above results. In particular if $R$ is a  reduced injective unstable module the localization $R \ra \ell(R)$ is an isomorphism, also if $I$ is injective analytic functor \cite{HLS93} then $f \circ m(I )\ra I$ is an equivalence.

\end{section}

\begin{section}{Proofs}

One reformulates \ref{uf} in a slightly different way than in the introduction.

\begin{thm}\label{uf}
Let $M$ be a finitely generated unstable module, there is  an  injective resolution $\mathcal I^\bullet$ of $M$ in $\U$ such that  any $\mathcal I^k$ is
a finite direct
sum of modules of the type $J(n) \otimes H^*V$.   
\end{thm}

In the introduction it is formulated  for a minimal injective resolution. In this setting the result will be that  $\mathcal I^k$ is a finite direct sum of indecomposable injective unstable $\A$-modules. Such modules are known to be of 
the form $J(n) \otimes L_\lambda$, $L_\lambda$ being an indecomposable factor of some $H^*V$, \cite{LS89} and see below. The two formulations are equivalent because $H^*V$ is a finite direct sum of indecomposable modules. The gain with the first one  is that  it allows to use the functor $T_V$, instead of the division functor by indecomposable factors $L_\lambda$ (see below). For what is necessary to the proofs  all division  functors share the same essential properties, however  $T_V$ is much more manageable.

The indecomposable injective unstable modules are, as is said above, of the form $J(n) \otimes L_\lambda$, where $L_\lambda$ is a direct factor in some $H^*V$. The isomorphism classes of indecomposable reduced injective unstable modules  are indexed by the simple representations over $\F2$ of the groups $\GL_n(\F2)$. These  representations are themselves indexed by
 $2$-regular partitions $\lambda$ of all integers,  $2$-regular partitions  being  the strictly decreasing ones. The unstable module $J(n) \otimes L_\lambda$ is the injective hull in $\U$ of of the unstable module $\Sigma^n S_\lambda(F(1))$, where $S_\lambda$ is the simple functor
in the category $\calF$ associated to $\lambda$. This functor is of degree $\vert \lambda \vert$ (\cite{PS98, Dj07, K94-2}).

\begin{thm}\label{lf} Let $M \in \U_n$ be finitely generated. Then the derivatives of $\ell$, for $i>0$,   $\ell^i(M)$ belong to $ \U_{n-1}$. Moreover the unstable modules $\ell^i(M)$ are finitely generated unstable modules. 
\end{thm}

The first part of the theorem is classical and follows from  \ref{krull} and the commutation of $T$ with $\ell$.

The next result is the stabilization theorem:

\begin{thm}\label{stab}
Let $M$ and $N$ be two finitely generated unstable modules. Then for any $i$ the $\F2$-vector space $\Ext^i_\calF(f(M),f(N))$ is  the colimit over $k$
of the system:
$$
\ldots \ra \ldots \Ext^i_\U(\Phi^kM,N) \ra \Ext^i_\U(\Phi^{k+1}M,N) \ra \ldots
$$
where the morphisms are induced by the maps $\Phi^{k+1}M \ra \Phi^{k}M$. 
Moreover for $k$ large enough there is an isomorphism:
$$
\Ext^i_\U(\Phi^kM,N)  \cong \Ext^i_\calF(f(M),f(N))
$$ \end{thm}

\vskip 1mm \vskip 1mm

The proofs of \ref{uf} and \ref{lf}  will be done at the same time by induction over the Krull filtration.
To prove \ref{uf} it is enough to prove:
\begin{prop} \label{uff} Let $M$ be a finitely generated unstable module and $k$ an integer. There exists only a finite number of unstable modules $\Sigma^n S_\lambda(F(1))$
such that $\Ext^k_\U(\Sigma^n S_\lambda(F(1)), M)$ is non trivial, and if it is non-trivial it is finite dimensional.
\end{prop}
Indeed,  $J(n) \otimes L_\lambda$ is the injective hull of $\Sigma^n S_\lambda(F(1))$, and  the dimension of the $\F2$-vector space $\Ext^k_\U(\Sigma^n S_\lambda(F(1)), M)$  is greater or equal to the number of occurrences of $J(n) \otimes L_\lambda$ in the term
$\mathcal I^k$ of a minimal resolution of $M$.

It is useful note, and one will use:
\begin{lem} \label{thick} If one is given a short exact sequence  $\{0\} \ra M' \ra M \ra M" \ra \{0\}$ if \ref{uff} and \ref{lf} hold for two terms of the sequence they hold for the third.
\end{lem}

\vskip 1mm \vskip 1mm

{\bf First step  of the proofs of \ref{uff} and \ref{lf}: $\U_0$.}

\begin{lem}
Both \ref{uff} and \ref{lf} are true for finitely generated locally finite unstable modules.
\end{lem}

Indeed, this follows from the following fact. Any finitely generated locally finite unstable modules has an injective resolution of finite length, such that any term of the resolution is a finite direct sum of Brown-Gitler modules (see \cite{S94}). More precisely one  constructs a finite resolution of a finitely generated object $M$ in
$\U_0$ because the injective hull $E_M$ of $M$ is also finite and thus finitely generated. Moreover, if one denotes by $v(M)$ the largest integer such that  $M$ is trivial in degrees strictly larger than $v(M)$, then
$v(E_M/M)<v(M)$. This allows to show that minimal resolution are finite.

For what concerns the localization functor and its derivatives they are all trivial, except if the module is concentrated in degree $0$, in which case the localization is an isomorphism, and the derivatives are trivial.

\vskip 1mm \vskip 1mm

{\bf Second step  of the proofs of \ref{uf} and \ref{lf}: $F(1)$.}

\begin{lem}\label{F(1)}
Both \ref{uf} and \ref{lf} are true for the unstable module  $F(1)$.
\end{lem}

This relies on computations in Mac Lane homology and explained in great details in \cite{NTC14}. We give  some informations below that are enough for our purpose.
 
 The unstable modules $\ell^i(F(1))$ are known to belong to $\U_0$ \ref{krull-loc}, it remains to prove they are finitely generated. In fact,  they are explicitly known, computations depending on Mac Lane homology show that
 they are finite and thus finitely generated. 
 This follows from theorem 12.13 of  "Alg\`ebre de Steenrod, modules instables et foncteurs polyn\^omiaux", \cite{ER03}. This allows to have some control on the $\mathcal N^k$ using \ref{4-terms} and to prove by
 induction they are finite. 
  It follows directly from these results that the nilpotent part of the minimal resolution of $F(1)$ is in each degree a finite direct sum of Brown-Gitler modules.
  
  It remains to show the result for the reduced part: that $\mathcal R^k$ of the  minimal resolution of $F(1)$ is a finite direct sum of indecomposable injective unstable modules. This is a consequence of the corresponding result in the category $\mathcal F$
\ref{pf}, \ref{puf}. 
  
  \vskip 2mm 
  
 On the way one gives some informations on the reduced part of the resolution of $F(1)$. One starts by a theorem which is proved in \cite{NTC14}. This is not necessary here but worth to be mentioned and the material introduced here, which is the longest part of the digression, will be used later. 
  
%
%

  \begin{prop}[A. Touz\'e {\cite[4.18]{To13}}]\label{To}
 Let $S$ be  a simple functor, if  $\Ext^*_\calF(S,Id)$ is non trivial then the degree of $S$ must be a power of $2$.
  \end{prop}
  
   One offers a proof different from the one of \cite{To13}, it depends on:
   
 \begin{lem}\label{spf}  Let $\lambda$ be a $2$-regular  partition, such that $\vert \lambda \vert =n$.  The functor $\Lambda^\lambda\cong\otimes  \Lambda^{\lambda_1} \otimes \ldots \Lambda^{\lambda_t}$ has  a filtration whose sub-quotients  are either simple functors of degree $m$ such that there exists $h$ with $2^hm=n$, or
   tensor products $F \otimes G$ of functors with no constant part ($F(\{0\})=G(\{0\})=\{0\}$).
   \end{lem}
To prove the proposition  consider  $\Lambda^\lambda \cong \Lambda^{\lambda_1}\otimes \ldots \otimes \Lambda^{\lambda_t}$ as a strict polynomial functor and one uses the classification of simple objects in $\mathcal P$. It follows from a theorem of Steinberg that they  are of the form, \cite{K02}
  $$
  S_{\lambda^{1}}^{(i_1)} \otimes   S_{\lambda^{1}}^{(i_2)} \otimes \ldots \otimes   S_{\lambda^{t}}^{(i_t)}
  $$ where $\lambda^{1}, \ldots , \lambda^{t}$ are $2$-regular partitions, $i_1 < i_2 < \ldots < i_t$, and  $S_{\lambda^{}}^{(i)}$ is the $i$-th Frobenius twist of the canonical lifting of 
  the simple functor $S_\lambda \in \calF$ to $\calP$.

 \vskip 1mm \vskip 1mm
 
 The following more precise form of \ref{spf} is necessary to complete the proof and will be necessary later: 
      \begin{lem}\label{basis} 
  Let  $\lambda$ be a $2$-regular partition of the integer $n$. The composition series of the strict polynomial functor $\Lambda^\lambda$ has
  \begin{itemize}
  \item  one sub-quotient $S_\lambda$,
  \item all other simple sub-quotients of degree $n$ as polynomial functors
  are functors $S_\mu$ with $\mu > \lambda$, for the natural order on partitions;
  \item Frobenius twists of functors $S_\nu^{(r)}$, with $2^r\vert \nu \vert =n$;
  \item non-trivial tensor products. 
 \end{itemize}
  \end{lem}
  
  This result can be deduced easily from various publications, in particular  \cite{PS98, Dj07} which do not claim for originality.
  \vskip 2mm
  
  Proposition \ref{To} follows by induction. As   $F$ and $G$ are functors with no constant part,  $\Ext^*_\calF(F \otimes G ,Id)$ is trivial. Hence, if $\Ext^*_\calF(\Lambda^{n},Id)$ is trivial, so is $\Ext^*_\calF(S_\lambda,Id)$    for all simple functors of degree $n$ by an increasing induction on the degree and a descending one  on the $2$-regular partitions of $n$. The result
  follows because $\Ext^*_\calF(\Lambda^{n},Id)$ is trivial if (and only if) $n$ is not a power of $2$ \cite{FLS94}.
  
   \vskip 2mm
   
   {\bf Third step  of the proofs of \ref{uff} and \ref{lf}: tensor products and $\U_1$.}
   
   From now on we  have proved the results for locally finite unstable modules and $F(1))$. To finish the case of $\U_1$ we will prove that if the theorems hold for two modules they hold for their tensor product. This follows from a standard double-complex spectral sequence argument. Then one uses \ref{thick} to reduce to this case (this uses a structural result from \cite{S98}).
   
   \begin{prop}\label{tp}
 If \ref{uff} and \ref{lf} hold  for $M$ and $N$,  they hold for the tensor product $M \otimes N$.
 \end{prop}
 
Theorem  \ref{lf} is true in this case,  because of  \ref{kunneth},  it is true for the tensor product $M \otimes N$. 

For theorem   \ref{uff},  let  $\mathcal I^\bullet$ and $\mathcal J^\bullet$ be injective resolutions of $M$ and $N$ having the required property.   The tensor product 
$\mathcal H^\bullet = \mathcal I^\bullet \otimes \mathcal J^\bullet$ is not an injective resolution of $M \otimes N$. Construct a Cartan-Eilenberg 
resolution $\mathcal H^{\bullet, \bullet}$ of $\mathcal H^\bullet$. Applying to this double complex the functor $\Hom_\U(S,-)$ yields an hypercohomology spectral sequence converging 
to $\Ext^{p+q}_\U(S, M \otimes N)$ with $E_2$-term $\Ext_\U^p(S, \mathcal H^q)$.

The unstable module $\mathcal H^{q}$ is a finite direct sum of modules $J(k) \otimes J(\ell) \otimes H^*V$. Thus, the group $\Ext^p_\U(S, \mathcal H^q) $ is isomorphic to a finite direct sum
of groups  $\Ext^p_\U(T_V(S),  J(k) \otimes J(\ell))$, the $\F2$-vector space $V$  being of bounded dimension as there are only finitely many factors. For $i=p+q$ let $d$ be an upper bound of the dimensions.

 Suppose now that $S$ is of the form $\Sigma^n S_\lambda(F(1))$. These groups are trivial
as soon as $n> sup(k+\ell)$ or as soon as the connectivity of  $T_V(S_\lambda(F(1))$,   $c(T_V(S_\lambda(F(1)))$ is greater than $  sup(k+\ell)$.  
 
 \begin{lem}The connectivity  $c(S_\lambda(F(1)))$,  of $S_\lambda(F(1))$, 
 is $\lambda'_1+2 \lambda'_2+ \ldots +2^{t-1}\lambda'_{t-1}-1 \geq \lambda'_1-1$, $\lambda'$ being the conigate (or dual) partition of $\lambda$.
 \end{lem}
 
 Using the properties of $T_V \cong T^{\dim(V)}$, the decomposition $T \cong \F2 \oplus \bar T$. and the fact that 
$S_\lambda (F(1)) \subset \tilde H^*(\Z/2)^{\lambda'_1}$, one shows that $c(T_V(S_\lambda(F(1))) \geq \lambda'_1-d-1$.

As a consequence,  for a given $V$ only 
a finite number of reduced simple unstable modules $T_VS_\lambda(F(1))$ are of connectivity less than 
a given constant. 

The result follows.

\vskip 1mm \vskip 1mm

This proves the theorems for $\U_1$. Indeed, they are  true for any $L \in\U_0$, $F(1)$ and by tensor product for any $L \otimes F(1)$. Next any finitely generated object  $M \in \U_1$ enters in a short exact sequence \cite{S98}, see also \cite{K13}:
$$
\{0\} \ra L' \ra M \ra L \otimes F(1) \ra L" \ra \{0\}
$$
where $L,L'q,L" \in \U_0$ and are finitely generated (see \cite{K13} for a generalization). The result follows.

\vskip 2mm

From now on one assumes that the theorems have been proved for objects in $\U_{n-1}$.

 \vskip 1mm \vskip 1mm
   
   {\bf Fourth  step  of the proofs of \ref{uff} and \ref{lf}: The case of $\Lambda^n(F(1))$.}

The following step is  the case of exterior powers. One assume the theorems  hold for $\Lambda^k(F(1))$, and prove  it holds for $\Lambda^{k+1}(F(1))$.

If $k$ is even $\Lambda^{k+1}(F(1))$ is  a direct summand in
$\Lambda^k(F(1)) \otimes F(1)$, and \ref{tp} implies that \ref{uff} and \ref{lf} hold for  $\Lambda^{k+1}(F(1))$. 
\vskip 1mm \vskip 1mm

If $k$ is odd the situation is more complicated. One  one has a short exact sequence
 $$\{0\} \ra W_{(k,1)}(F(1)) \ra \Lambda^k(F(1)) \otimes F(1) \ra \Lambda^{k+1}(F(1)) \ra \{0\}$$
 which defines $\{0\} \ra W_{(k,1)}(F(1))$. There is another short exact sequence:
 $$\{0\}\ra \Lambda^{k+1}(F(1)) \ra W_{(k,1)}(F(1)) \ra S_{(k,1)}(F(1))\ra \{0\} $$
   which defines $\{0\} \ra S_{(k,1)}(F(1))$. 
By dualization (\cite{PS98}, \cite{Dj07}) one gets
 $$\{0\} \ra \Lambda^{k+1}(F(1)) \ra \Lambda^k(F(1)) \otimes F(1) \ra  DW_{(k,1)}(F(1)) \ra \{0\}$$
 $$\{0\}\ra S_{(k,1)}(F(1))\ra DW_{(k,1)}(F(1)) \ra \Lambda^{k+1}(F(1))  \ra \{0\} $$. Recall (see the above mentioned references that the simple functor $S_{(k,1)}$ is self-dual.
 
 These exact sequences define $ W_{(k,1)}(F(1))$ and $ DW_{(k,1)}(F(1))$,  which defines $S_{(k,1)}(F(1))$. .

Assume (for the integer $k$) that $\ell^i( \Lambda^{k+1}(F(1)))$, $\ell^i( S_{(k,1)}(F(1)))$ are finitely generated and show the result for $k+1$.
The same is true for $\ell^i(W_{(k,1)}(F(1)))$ and $\ell^i(DW_{(k,1)}(F(1)))$.

Using  the long exact sequences associated to the above   exact sequences one first shows that 
$\ell^{i}(W_{(k+1,1)}(F(1)))$ is finitely generated. Then one gets the result for  $\ell^{i}(\Lambda^{(k+2)}(F(1)))$, next for $\ell^{i}(S_{(k+1,1)}(F(1)))$ and finally for 
$\ell^{i}(DW_{(k+1,1)}(F(1)))$. This last case is not used in the induction.

\vskip 1mm \vskip 1mm

The case of the exterior powers follows also from explicit  computations  \cite{FFSS99}.

\vskip 1mm \vskip 1mm

 {\bf Last  step  of the proofs of \ref{uff} and \ref{lf}.}

 As the result holds for all exterior powers $\Lambda^k(F(1))$, it holds
for tensor products of such modules. Then,  \ref{basis} implies it holds for any  module $S_\lambda(F(1))$ 
$\vert \lambda \vert=n$ and using 
\ref{tp} for all $\Sigma^nS_\lambda(F(1))$.

To finish the proof  one uses the following two lemmas and the induction hypothesis on the Krull filtration. 

\begin{lem}\label{Krull-Omega}
If a reduced unstable module $M$ belongs to $\U_n$, $\Omega M$ belongs to $\U_{n-1}$.
\end{lem}

The proof is left to the reader.

\begin{lem}
If the theorem holds for a reduced unstable module $M$ of Krull filtration $n$, and for all unstable modules of Krull fitration $n-1$, then it holds
for any submodule $N$ of $M$ so that there
 exists $\ell$ with $\Phi^\ell M \subset N \subset M$.
\end{lem}

This is because $M/\Phi^\ell M \in \U_{n-1}$.

\vskip 1mm \vskip 1mm

The proof ends  using \ref{KS} and the following:

\begin{prop} A finitely generated reduced unstable $R$ of Krull filtration $n$ has a finite filtration whose quotients are suspensions of reduced modules of Krull filtration less than $n$, and whose
associated functors are simple (of degree less than $n$). 
\end{prop}
This means that the sub-quotients are of the form
$\Sigma^n R$, for some $R$ so that $\Phi^kS_\lambda \subset R \subset S_\lambda$ for some $\lambda$.  The required filtration is obtained as follows. One considers the localization $R \ra \ell(R)$. Following Kuhn, \cite{K94-2}, one knows that  $f(R)=f(\ell(R))$ has a finite composition series. The proof is done by induction on the length of the composition series of $f(R)$. Consider an epimorphism on a  simple functor: $f(R) \ra S_\lambda$, and the associated
unstable module map $R \ra S_\lambda(F(1))$. Its kernel $K$  is reduced, and  one can apply the induction hypothesis.

\vskip 1mm
To finish the proof of the theorem one applies the preceding lemma.

\vskip 1mm

{\bf Proof of \ref{stab} and of the corollary}

The proof of \ref{stab} follows directly from \ref{uff} and from    the properties of the right adjoint $\tilde \Phi$ of $\Phi$, and in particular from the computation
$\tilde \Phi(L \otimes J(n) \cong L \otimes J(n/2)$, ($L$ a reduced injective, and $J(n/2)$ is trivial if $n$ is not an integer). If ${\mathcal N}^i$ is a finite direct sum of indecomposable injective unstable modules it is clear that
for $k$ large enough $\tilde \Phi^k ({\mathcal N}^i)$ is trivial. The use of $\tilde \Phi$ can be replaced by:

\begin{prop}
Let $M$ be a finitely generated unstable modules, assume $N$ is nilpotent and has finite nilpotent filtration. If $k$ is large enough $\Hom_\U(\Phi^k(M),N) \cong\{0\}$.
\end{prop}

Note that in particular a finitely generated module has a finite nilpotent filtration \cite{S94}, \cite{K95}. Suppose $N_t \{0\} \cong \{0\}$ for $t>t_0$.  The unstable module $\Omega \Phi^k(M)$ 
is isomorphic to $\Sigma^{2^k-1}\Phi^k \Omega M$, and has trivial image in $N$ if $2^k-1 >t_0$.

For the corollary, one observes that \ref{stab} allows to reduce to a computation in the category $\calF$. It is enough to show the following:

\begin{prop}
Let $F$ be a polynomial functor of degree $n$, such that $F(\{0\})=\{0\}$, $R$ the smallest sub-functor such that $F/R$ is of degree $n-1$. 
Then $\Ext_\calF^1( R \otimes R , F) \cong \{0\}$.
\end{prop}

This last result is proved using the category $bi-\calF$, \cite{FFSS99}. This result is used in\cite{CGS14}.

\end{section}

\bibliographystyle{amsalpha}

\providecommand{\bysame}{\leavevmode\hbox to3em{\hrulefill}Thinspace}
\providecommand{\MR}{\relax\ifhmode\unskip\space\fi MR }
\providecommand{\MRhref}[2]{%
  \href{http://www.ams.org/mathscinet-getitem?mr=#1}{#2}
}
\providecommand{\href}[2]{#2}

\end{document}